\documentclass[11pt]{article}

\usepackage[a4paper,margin=31mm]{geometry}
\usepackage{amsmath,amssymb,amsthm,mathtools}
\usepackage[T1]{fontenc}
\usepackage{lmodern}
\usepackage{microtype}
\usepackage{enumitem}
\usepackage[hidelinks]{hyperref}

\newtheorem{theorem}{Theorem}[section]
\newtheorem{proposition}[theorem]{Proposition}
\newtheorem{lemma}[theorem]{Lemma}
\newtheorem{corollary}[theorem]{Corollary}
\newtheorem{remark}[theorem]{Remark}

\newcommand{\B}{\mathcal B}
\newcommand{\K}{\mathcal{K}}
\newcommand{\Hq}{\mathcal{H}_0}
\newcommand{\Ran}{\mathcal{R}}
\newcommand{\Ker}{\mathcal{N}}
\newcommand{\PF}{\operatorname{PF}}
\newcommand{\LGD}{\operatorname{LGD}}
\newcommand{\RGD}{\operatorname{RGD}}
\newcommand{\GD}{\operatorname{GD}}
\newcommand{\Red}{\operatorname{Red}}
\newcommand{\spec}{\sigma}

\title{Pseudo-Fredholm operators form a regularity on Hilbert spaces}

\author{Zakariae Aznay\textsuperscript{1}\quad
Fadil Chabbabi\textsuperscript{1}\quad
Hassan Zariouh\textsuperscript{2}\\[8pt]
\normalsize\textsuperscript{1}Laboratory (A.F.N.L.A.), Department of Mathematics,\\[-1pt]
\normalsize Faculty of Sciences, Abdelmalek Essa\^adi University,\\[-1pt]
\normalsize Tetouan, Morocco\\[2pt]
\normalsize E-mail (Z. Aznay): \texttt{aznay.zakariae@gmail.com}\\[-1pt]
\normalsize E-mail (F. Chabbabi): \texttt{f.chabbabi@uae.ac.ma}\\[7pt]
\normalsize\textsuperscript{2}Department of Mathematics (CRMEFO), Laboratory (L. CS2A),\\[-1pt]
\normalsize Faculty of Science, Mohammed I University,\\[-1pt]
\normalsize 60000 Oujda, Morocco\\[2pt]
\normalsize E-mail (H. Zariouh): \texttt{h.zariouh@yahoo.fr}}

\date{}

\begin{document}

\maketitle

\begin{abstract}
Let $H$ be a complex Hilbert space. We prove that the class of
pseudo-Fredholm operators on $H$ is a regularity in the sense of
Kordula and M\"uller. We also prove that the classes of left and right
generalized Drazin invertible operators are regularities on $H$.
Consequently, the generalized Kato spectrum and the corresponding
generalized Drazin spectra satisfy the spectral mapping theorem.
We construct an operator $T$ on a separable complex Banach space
such that $T^2$ is right generalized Drazin invertible whereas $T$
is not pseudo-Fredholm. In particular, the power axiom need not hold
for pseudo-Fredholm operators on Banach spaces. Finally, we show
that the accumulation sets of the approximate point and surjectivity
spectra define regularities on Banach spaces, yielding another proof
of the regularity of generalized Drazin invertible operators in $\B(X)$.
\end{abstract}


\section{Introduction and preliminaries}

Let $X$ be a complex Banach space and let $\B(X)$ denote the
algebra of bounded linear operators on $X$. For $T\in\B(X)$,
we denote by $\Ker(T)$ and $\Ran(T)$ its kernel and range,
respectively. For $A\in\B(X)$ and $M\subseteq X$, we write
$A^{-1}(M)=\{x\in X:Ax\in M\}$ for the inverse image of $M$ under $A$.
A pair $(M,N)$ belongs to $\Red(T)$ if $M$ and
$N$ are closed $T$-invariant subspaces and $X=M\oplus N$.
For a closed $T$-invariant subspace $M$, the notation $T|_M$
denotes the restriction of $T$ to $M$, viewed as an operator
from $M$ into $M$. All direct sums of closed subspaces
considered below are topological; they are not necessarily
orthogonal.

An operator $S\in\B(X)$ is called \emph{semi-regular} if
$\Ran(S)$ is closed and
\[
   \Ker(S)\subseteq
   \Ran^\infty(S):=\bigcap_{n\geq1}\Ran(S^n).
\]
According to Mbekhta~\cite{Mbekhta1990}, $T$ is
\emph{pseudo-Fredholm} if there exists $(M,N)\in\Red(T)$
such that $T|_M$ is semi-regular and $T|_N$ is
quasinilpotent; such a pair is called a generalized Kato
decomposition of $T$. We write $\PF(X)$ for the class of
pseudo-Fredholm operators acting on $X$.

Left and right generalized Drazin invertibility were introduced
in \cite{HocineBenharratMessirdi2015}. We use the following
equivalent descriptions: an operator is left (resp,
right) generalized Drazin invertible if it admits a reducing
decomposition into a bounded-below (resp, surjective)
operator and a quasinilpotent operator; see
\cite[Propositions~2.3 and~2.4]{BenharratHocineMessirdi2019}.
An operator is generalized Drazin invertible if it admits a
reducing decomposition into an invertible operator and a
quasinilpotent operator. We denote these three classes by
$\LGD(X)$, $\RGD(X)$ and $\GD(X)$, respectively.

Kordula and M\"uller~\cite[Definition~1.2]{KM1996}
introduced the notion of regularity in a complex unital
Banach algebra. In the present setting a nonempty subset
$\mathcal R\subseteq\B(X)$ is a regularity if
\begin{equation}\label{eq:regularity-powers}
 T\in\mathcal R
 \quad\Longleftrightarrow\quad
 T^n\in\mathcal R
 \qquad(n\geq1),
\end{equation}
and
\[
 AB\in\mathcal R
 \quad\Longleftrightarrow\quad
 A,B\in\mathcal R,
\]
whenever $A,B,C,D\in\B(X)$ commute pairwise and satisfy
$AC+BD=I$.

Lubansky~\cite{Lubansky2007} proved that the generalized
Drazin invertible elements of a complex unital Banach algebra
form a regularity. In particular, $\GD(X)$ is a regularity
for every complex Banach space $X$. We prove that $\PF(H)$
is a regularity when $H$ is a complex Hilbert space, and we 
construct a Banach space example for which the power axiom fails.

For $T\in\B(X)$, its analytic core and quasinilpotent
part are defined by
\begin{align*}
 \K(T)&=
 \left\{x\in X:
 \begin{array}{l}
 \text{there exist }c>0\text{ and }
 (x_n)_{n\geq0}\subseteq X\text{ such that}\\[-2pt]
 x_0=x,\quad Tx_{n+1}=x_n,\quad
 \|x_n\|\leq c^n\|x\|\quad(n\geq0)
 \end{array}
 \right\},\\
 \Hq(T)&=
 \left\{x\in X:
 \lim_{n\to\infty}\|T^nx\|^{1/n}=0
 \right\}.
\end{align*}
We use the standard properties
\begin{equation}\label{eq:core-facts}
 T\K(T)=\K(T),\qquad
 \K(T)\subseteq\Ran^\infty(T),\qquad
 T^{-1}\Hq(T)=\Hq(T).
\end{equation}
If $S$ is a semi-regular operator on a complex Hilbert space, then
$\K(S)=\Ran^\infty(S)$ is closed and
$\Hq(S)\subseteq\K(S)$, see
Mbekhta~\cite[Proposition~1.4(e), (g)]{Mbekhta1990}.
Moreover, $\K(Q)=\{0\}$ and $\Hq(Q)=X$ whenever
$Q$ is quasinilpotent.

We use the following characterization from
\cite[Corollary~3.2]{AznayChabbabi2026}: for every
complex Hilbert space $H$,
\begin{equation}\label{eq:hilbert-criterion}
 T\in\PF(H)
 \quad\Longleftrightarrow\quad
 \K(T),\quad
 \K(T)+\Hq(T),\quad
 \Ran(T)+\Hq(T)
 \ \text{are closed in }H.
\end{equation}
We give an alternative proof of this criterion by triangular
reduction. The proof of regularity uses this criterion and a
quotient argument involving left generalized Drazin invertibility.
The classes $\LGD(H)$ and $\RGD(H)$ are also shown to be regularities.
An earlier paper on spectral mapping for the generalized Kato
spectrum was retracted; see the retraction notice
\cite{Retraction2016}. No result from that paper is used here.
Our spectral mapping assertions follow from the regularity axioms.

The Banach-space counterexample below also concerns isolated points
of the surjectivity spectrum. We finally consider regularities
associated with the accumulation sets of one-sided spectra.

\section{Triangular reduction and quotients}

\begin{lemma}\label{lem:triangular}
Let $H$ be a complex Hilbert space and
$T\in\B(H)$. Then $T$ is pseudo-Fredholm if
and only if there exist closed $T$-invariant
subspaces $E\subseteq F$ such that $T|_E$
is surjective, the operator induced by $T$
on $F/E$ is quasinilpotent, and the operator
induced by $T$ on $H/F$ is bounded below.
\end{lemma}

\begin{proof}
Suppose first that such subspaces $E\subseteq F$
exist. Put $G=F\cap E^\perp$.
With respect to the orthogonal decomposition
$H=E\oplus G\oplus F^\perp$, write
\[
 T=
 \begin{pmatrix}
 S&U&W\\
 0&Q&V\\
 0&0&D
 \end{pmatrix}.
\]
Here $S=T|_E$ is surjective. Since $E$ and
$F$ are $T$-invariant, $Q=P_GT|_G$ and
$D=P_{F^\perp}T|_{F^\perp}$ are unitarily
equivalent to the operators induced by $T$
on $F/E$ and $H/F$, respectively, via the
canonical isometries $g\mapsto g+E$ and
$z\mapsto z+F$. Hence $Q$ is quasinilpotent
and $D$ is bounded below.

Since the spaces are Hilbert, there are bounded
operators $R,J$ such that
\[
 SR=I_E,\qquad JD=I_{F^\perp}.
\]
The series
\[
 X=-\sum_{n=0}^{\infty}R^{n+1}UQ^n,
 \qquad
 Y=\sum_{n=0}^{\infty}Q^nVJ^{n+1}
\]
converge absolutely in operator norm.
Indeed, since $Q$ is quasinilpotent, for
every $\varepsilon>0$ there exists
$M_\varepsilon>0$ such that
\[
 \|Q^n\|\leq M_\varepsilon\varepsilon^n
 \qquad(n\geq0).
\]
Choosing $\varepsilon$ so that
\[
 \varepsilon\max\{\|R\|,\|J\|\}<1
\]
proves the asserted convergence. Direct
calculation gives
\[
 SX-XQ=-U,\qquad QY-YD=-V.
\]
Put
\[
 P_1=
 \begin{pmatrix}
 I&X&0\\
 0&I&0\\
 0&0&I
 \end{pmatrix},
 \qquad
 P_2=
 \begin{pmatrix}
 I&0&0\\
 0&I&Y\\
 0&0&I
 \end{pmatrix}.
\]
Both $P_1$ and $P_2$ are bounded and invertible, and the
preceding identities give
\[
 P_2^{-1}P_1^{-1}TP_1P_2
 =
 \begin{pmatrix}
 S&0&W-XV\\
 0&Q&0\\
 0&0&D
 \end{pmatrix}.
\]
Hence $T$ is similar to $Z\oplus Q$, where
\[
 Z=
 \begin{pmatrix}
 S&W-XV\\
 0&D
 \end{pmatrix}
 \quad\text{ acts on }E\oplus F^\perp.
\]
Since $D$ is injective and $S$ is surjective,
\[
 \Ker(Z)=\Ker(S)\oplus\{0\}.
\]
Moreover, for any $u\in E$ and
$v=Dz\in\Ran(D)$, choose $e\in E$ such
that $Se=u-(W-XV)z$. Then
$Z(e,z)=(u,v)$, so that
\[
 \Ran(Z)=E\oplus\Ran(D),
\]
which is closed since $D$ is bounded below.
For every $n\geq1$,
\[
 Z^n(e,0)=(S^ne,0);
\]
the surjectivity of $S^n$ therefore yields
\[
 \Ker(Z)\subseteq E\oplus\{0\}
 \subseteq\bigcap_{n\geq1}\Ran(Z^n).
\]
Thus $Z$ is semi-regular. Since $Q$ is
quasinilpotent, $T$ is pseudo-Fredholm.

Conversely, suppose that $T=S\oplus Q$
is a generalized Kato decomposition on
a topological direct sum $H=M\oplus N$,
where $S$ is semi-regular and $Q$ is
quasinilpotent. Put $K=\K(S)$ and choose
\[
 E=K\oplus\{0\},\qquad F=K\oplus N.
\]
Since $S$ is semi-regular, $\K(S)$ is closed
and $S\K(S)=\K(S)$. Hence, as $H=M\oplus N$
is a topological direct sum, $E$ and $F$
are closed $T$-invariant subspaces, and
$T|_E$ is surjective. Since $F/E$ is
naturally isomorphic to $N$, the operator
induced by $T$ on $F/E$ is similar to $Q$
and hence quasinilpotent.

Furthermore,
\[
 \Ker(S)\subseteq K,\qquad S(K)=K
\]
imply $S^{-1}(K)=K$. Thus the operator
induced by $S$ on $M/K$ is injective,
with range $\Ran(S)/K$, which is closed
since $\Ran(S)$ and $K$ are closed and
$K\subseteq\Ran(S)$. This induced operator
is therefore bounded below. Since $H/F$
is naturally isomorphic to $M/K$, the
operator induced by $T$ on $H/F$ is also
bounded below.
\end{proof}

\begin{proposition}\label{prop:criterion}
For a complex Hilbert space $H$ and
$T\in\B(H)$, the three spaces
$\K(T)$, $\K(T)+\Hq(T)$ and
$\Ran(T)+\Hq(T)$ are closed if and only
if $T$ is pseudo-Fredholm.
\end{proposition}

\begin{proof}
If $T=S\oplus Q$ is a generalized Kato
decomposition, the facts recalled in
Section~1 give
\[
 \K(T)=\K(S)\oplus\{0\},\quad
 \K(T)+\Hq(T)=\K(S)\oplus N,\quad
 \Ran(T)+\Hq(T)=\Ran(S)\oplus N,
\]
where $S$ acts on $M$, $Q$ acts on $N$,
and $\K(S)$ and $\Ran(S)$ are closed.
Thus the three spaces are closed.

Conversely, suppose they are closed and put
\[
 E=\K(T),\qquad F=\K(T)+\Hq(T).
\]
Then $T(E)=E$. Every coset in $F/E$
has a representative $h\in\Hq(T)$;
therefore, if $Q$ denotes the induced
operator on $F/E$,
\[
 \|Q^n(h+E)\|\leq\|T^nh\|.
\]
Thus $\Hq(Q)=F/E$. Since $F/E$ is a Hilbert space,
$Q$ is quasinilpotent by \cite[Proposition~1.4(d)]{Mbekhta1990}.

Moreover, $T^{-1}(F)=F$: if $Tx=k+h$
with $k\in\K(T)$ and $h\in\Hq(T)$,
choose $u\in\K(T)$ with $Tu=k$
and use $T^{-1}(\Hq(T))=\Hq(T)$.
Hence the induced operator on $H/F$
is injective, with closed range
\[
 (\Ran(T)+\Hq(T))/F.
\]
It is therefore bounded below.
Lemma~\ref{lem:triangular} gives that 
$T\in\PF(H)$.
\end{proof}

\begin{lemma}\label{lem:left-criterion}
Let $H$ be a complex Hilbert space and
$T\in\B(H)$. Then
\[
 T\in\LGD(H)
 \quad\Longleftrightarrow\quad
 \Hq(T)\text{ and }\Ran(T)+\Hq(T)
 \text{ are closed}.
\]
\end{lemma}

\begin{proof}
If $T=D\oplus Q$, with $D$ bounded
below and $Q$ quasinilpotent, then
$\Hq(T)$ is precisely the subspace
on which $Q$ acts, while
\[
 \Ran(T)+\Hq(T)
 =\Ran(D)\oplus\Hq(T).
\]
Both spaces are closed.

Conversely, put $N=\Hq(T)$ and suppose
that $N$ and $\Ran(T)+N$ are closed.
Since $\Hq(T|_N)=N$, the restriction $T|_N$ is quasinilpotent
by \cite[Proposition~1.4(d)]{Mbekhta1990}.
By assertion \eqref{eq:core-facts} above, the operator
$\widehat T$ induced by $T$ on $H/N$
is injective, and
\[
 \Ran(\widehat T)
 =(\Ran(T)+N)/N
\]
is closed. Hence $\widehat T$ is bounded
below. Applying the proof of
Lemma~\ref{lem:triangular} with
$E=\{0\}$ and $F=N$ we obtain  a decomposition
of $T$ into a bounded-below operator
and a quasinilpotent operator.
\end{proof}

We shall also use the following quotient criterion.

\begin{lemma}\label{lem:quotient}
Let $H$ be a complex Hilbert space,
$T\in\B(H)$ and $E$ a closed subspace
such that $T(E)=E$. Then $T$ is
pseudo-Fredholm whenever the operator
induced by $T$ on $H/E$ is left
generalized Drazin invertible.

Conversely, if $T$ is pseudo-Fredholm,
then $E=\K(T)$ is closed, $T(E)=E$,
and the operator induced by $T$ on
$H/E$ is left generalized Drazin
invertible.
\end{lemma}

\begin{proof}
Suppose first that the induced operator
on $H/E$ is left generalized Drazin
invertible. Write
\[
 \widehat T=A\oplus Q
\]
on $H/E=M_0\oplus N_0$, where $A$
is bounded below and $Q$ is
quasinilpotent. Put $F=\pi^{-1}(N_0)$,
where $\pi:H\to H/E$ is the canonical
quotient map. Then $E\subseteq F$ are
closed $T$-invariant subspaces,
$T|_E$ is surjective, and the operators
induced by $T$ on $F/E$ and $H/F$
are similar to $Q$ and $A$,
respectively. The conclusion follows
from Lemma~\ref{lem:triangular}.

Conversely, suppose that $T=S\oplus Q$
on $H=M\oplus N$ is a generalized
Kato decomposition, where $S$ is
semi-regular and $Q$ is quasinilpotent.
Then
\[
 \K(T)=\K(S)\oplus\{0\}
\]
is closed and $T(\K(T))=\K(T)$.
As shown in the proof of Lemma~\ref{lem:triangular},
the operator $\widehat S$ induced by $S$ on $M/\K(S)$
is bounded below.
Finally, the natural identification
\[
 H/\K(T)\simeq (M/\K(S))\oplus N
\]
identifies the operator induced by $T$
with the left generalized Drazin
invertible operator $\widehat S\oplus Q$.
\end{proof}

\section{Regularity of pseudo-Fredholm operators}

The next lemma is used to decompose the quasinilpotent part
of a product in Lemma~\ref{lem:bezout}.

\begin{lemma}\label{lem:local-quasinilpotent}
Let $X$ be a complex Banach space, $T\in\B(X)$,
and $x\in \Hq(T)$. Then there exists a linear subspace
$Y_x$ of $X$ containing $x$, equipped with a norm
$\|\cdot\|_x$, such that:
\begin{enumerate}
\item[(i)] $(Y_x,\|\cdot\|_x)$ is a Banach space
continuously embedded in $X$, with
\[
\|y\|\leq\|y\|_x\qquad(y\in Y_x).
\]
\item[(ii)] every $R\in\B(X)$ commuting with $T$
leaves $Y_x$ invariant and satisfies
\[
\|R|_{Y_x}\|_{\B(Y_x)}\leq\|R\|.
\]
\item[(iii)] $T|_{Y_x}$ is quasinilpotent.
\end{enumerate}
In particular, $x\in Y_x\subseteq \Hq(T)$.
\end{lemma}
\begin{proof}
For $n\geq1$, put
\[
 r_n=\frac1n+\sup_{j\geq n}\|T^jx\|^{1/j},
 \qquad
 w_0=1,\qquad w_n=r_1\cdots r_n.
\]
Since $x\in\Hq(T)$, the sequence $(r_n)$ is positive,
decreasing and converges to zero. Define
\[
 Y_x=\left\{y\in X:
 \sup_{n\geq0}\frac{\|T^ny\|}{w_n}<\infty\right\},
 \qquad
 \|y\|_x=\sup_{n\geq0}\frac{\|T^ny\|}{w_n}.
\]
For $1\leq j\leq n$, we have
$r_j\geq\|T^nx\|^{1/n}$, so $\|T^nx\|\leq w_n$.
Consequently,
\[
 \|x\|_x\leq\max\{\|x\|,1\}<\infty,
\]
and hence $x\in Y_x$.

For $y,z\in Y_x$ and $\alpha,\beta\in\mathbb C$,
\[
 \sup_{n\geq0}
 \frac{\|T^n(\alpha y+\beta z)\|}{w_n}
 \leq|\alpha|\|y\|_x+|\beta|\|z\|_x.
\]
Thus $Y_x$ is a linear subspace of $X$.
The homogeneity of $\|\cdot\|_x$ is immediate, and
the preceding inequality gives the triangle inequality.
Moreover, $w_0=1$ implies
\[
 \|y\|\leq\|y\|_x\qquad(y\in Y_x).
\]
Therefore $\|\cdot\|_x$ is a norm and the inclusion
$Y_x\hookrightarrow X$ is continuous.

Let $(y_m)$ be a Cauchy sequence in $Y_x$. Since
$\|y_m-y_\ell\|\leq\|y_m-y_\ell\|_x$, it is  a Cauchy sequence in $X$.
The completeness of $X$ gives $y\in X$ such that $y_m\to y$
in $X$. Given $\varepsilon>0$, choose
$m_0$ such that
\[
 \|T^n(y_m-y_\ell)\|\leq\varepsilon w_n
 \qquad(m,\ell\geq m_0,\ n\geq0).
\]
For fixed $m\geq m_0$ and $n\geq0$, letting
$\ell\to\infty$ gives
\[
 \|T^n(y_m-y)\|\leq\varepsilon w_n.
\]
It follows that
\[
 \sup_{n\geq0}\frac{\|T^ny\|}{w_n}
 \leq\|y_{m_0}\|_x+\varepsilon<\infty.
\]
Hence $y\in Y_x$ and $\|y_m-y\|_x\leq\varepsilon$
for $m\geq m_0$. Thus $Y_x$ is complete.

If $R\in\B(X)$ commutes with $T$, then
\[
 \sup_{n\geq0}\frac{\|T^nRy\|}{w_n}
 \leq\|R\|\sup_{n\geq0}\frac{\|T^ny\|}{w_n}
 =\|R\|\,\|y\|_x.
\]
Therefore $R(Y_x)\subseteq Y_x$ and
$\|R|_{Y_x}\|_{\B(Y_x)}\leq\|R\|$.

Finally, the monotonicity of $(r_n)$ gives
\[
 \frac{w_{n+k}}{w_n}
 =r_{n+1}\cdots r_{n+k}
 \leq r_1\cdots r_k=w_k
 \qquad(n\geq0,\ k\geq1).
\]
Consequently,
\[
 \|T^ky\|_x
 \leq
 \left(\sup_{n\geq0}\frac{w_{n+k}}{w_n}\right)\|y\|_x
 \leq w_k\|y\|_x.
\]
Since $r_n\to0$, we have $w_k^{1/k}\to0$. Hence
\[
 \|(T|_{Y_x})^k\|_{\B(Y_x)}^{1/k}
 \leq w_k^{1/k}\longrightarrow0,
\]
which proves that $T|_{Y_x}$ is quasinilpotent.
Moreover, for every $y\in Y_x$,
\[
 \|T^ny\|\leq w_n\|y\|_x.
\]
Taking $n$th roots shows that $y\in\Hq(T)$.
\end{proof}

\begin{lemma}\label{lem:bezout}
Let $X$ be a complex Banach space and let
$A,B,C,D\in\B(X)$ commute pairwise with $AC+BD=I$.
Then
\begin{align}
 \K(AB)&=\K(A)\cap\K(B),
 \label{eq:bezout-K}\\
 \Hq(AB)&=\Hq(A)\oplus\Hq(B),
 \label{eq:bezout-H}\\
 \Hq(A)&\subseteq\K(B),\qquad
 \Hq(B)\subseteq\K(A),
 \label{eq:bezout-cross}\\
 \Ran(AB)&=\Ran(A)\cap\Ran(B).
 \label{eq:bezout-range}
\end{align}
The sum in~\eqref{eq:bezout-H} is algebraic.
If $\Hq(AB)$ is closed, then its two summands are
closed and the sum is topological. Moreover, putting
\[
 L(U)=\K(U)+\Hq(U),\qquad
 G(U)=\Ran(U)+\Hq(U),
\]
we have
\[
 L(AB)=L(A)\cap L(B),\qquad
 G(AB)=G(A)\cap G(B),
\]
and
\begin{equation}\label{eq:bezout-G-preimages}
 A^{-1}G(AB)=G(B),\qquad
 B^{-1}G(AB)=G(A).
\end{equation}
\end{lemma}

\begin{proof}
Identity~\eqref{eq:bezout-K} is the subject of 
\cite[Lemma~2.3(i)]{Tajmouati2014}.
We prove the remaining assertions.

Put $T=AB$, $U=AC$ and $V=BD$. Then
\[
 U+V=I,\qquad UV=TCD.
\]
Let $x\in\Hq(T)$. By
Lemma~\ref{lem:local-quasinilpotent}, there exists
a Banach space $Y_x$, continuously embedded in $X$
and containing $x$, such that $T|_{Y_x}$ is
quasinilpotent and the restrictions of $A,B,C,D$
to $Y_x$ are bounded. In the following spectral argument,
$Y_x$ carries the norm $\|\cdot\|_x$, and its subspaces carry
the induced norm. Since
\[
 (U|_{Y_x})(I_{Y_x}-U|_{Y_x})
 =(T|_{Y_x})((CD)|_{Y_x})
\]
is a product of commuting operators, one of which is quasinilpotent,
it is quasinilpotent. Polynomial spectral mapping gives
\[
 \spec(U|_{Y_x})\subseteq\{0,1\}.
\]
The Riesz decomposition theorem yields
\[
 Y_x=Y_{x,0}\oplus Y_{x,1},
\]
where
\[
 \spec(U|_{Y_{x,0}})\subseteq\{0\},
 \qquad
 \spec(U|_{Y_{x,1}})\subseteq\{1\},
\]
and either summand may be zero. These subspaces are
closed in $Y_x$ (not necessarily in $X$) and invariant under $A,B,C,D$, since
their restrictions commute with the spectral
projections of $U|_{Y_x}$.

On $Y_{x,0}$, the operator
$V|_{Y_{x,0}}=I_{Y_{x,0}}-U|_{Y_{x,0}}$ is invertible.
Since $BD=DB=V$, the restriction $B|_{Y_{x,0}}$
is invertible, with inverse
$(D|_{Y_{x,0}})(V|_{Y_{x,0}})^{-1}$. Thus
\[
 A|_{Y_{x,0}}
 =(T|_{Y_{x,0}})(B|_{Y_{x,0}})^{-1}
\]
is quasinilpotent. Similarly, $A|_{Y_{x,1}}$ is
invertible and $B|_{Y_{x,1}}$ is quasinilpotent.
For $i=0,1$, let $\|\cdot\|_{x,i}$ denote the operator norm
on $Y_{x,i}$. For $y\in Y_{x,0}$,
\[
 \|A^ny\|\leq\|A^ny\|_x
 \leq\|(A|_{Y_{x,0}})^n\|_{x,0}\,\|y\|_x.
\]
Taking $n$th roots, and arguing similarly on $Y_{x,1}$, gives
\[
 Y_{x,0}\subseteq\Hq(A),\qquad
 Y_{x,1}\subseteq\Hq(B).
\]
Decomposing $x=x_0+x_1$ with $x_i\in Y_{x,i}$,
we obtain $x\in\Hq(A)+\Hq(B)$. The reverse inclusion
follows from commutativity.

To prove directness, let $x\in\Hq(A)\cap\Hq(B)$
and set $M=\max\{1,\|U\|,\|V\|\}$. Then
\[
 \|U^nx\|^{1/n}\longrightarrow0,\qquad
 \|V^nx\|^{1/n}\longrightarrow0.
\]
For $n\leq k\leq2n$,
\[
 \|U^kV^{2n-k}x\|
 =\|U^{k-n}V^{2n-k}U^nx\|
 \leq M^n\|U^nx\|,
\]
whereas for $0\leq k<n$,
\[
 \|U^kV^{2n-k}x\|
 =\|U^kV^{n-k}V^nx\|
 \leq M^n\|V^nx\|.
\]
Applying the binomial formula and using $U+V=I$,
we obtain
\[
 \|x\|=\|(U+V)^{2n}x\|
 \leq(4M)^n\bigl(\|U^nx\|+\|V^nx\|\bigr)
 \longrightarrow0.
\]
Thus $x=0$, proving~\eqref{eq:bezout-H}.

Suppose that $E=\Hq(T)$ is closed. Then $T|_E$ is
quasinilpotent. Applying the preceding spectral argument
on $E$ gives a topological decomposition
\[
 E=E_0\oplus E_1,\qquad
 E_0\subseteq\Hq(A),\quad E_1\subseteq\Hq(B).
\]
By~\eqref{eq:bezout-H} and uniqueness of the
decomposition, $E_0=\Hq(A)$ and $E_1=\Hq(B)$.
Hence both spaces are closed in $X$.

Next, let $x\in\Hq(A)$. Applying
Lemma~\ref{lem:local-quasinilpotent} to $A$, we obtain
a Banach space $Y_x$ on which $A$ is quasinilpotent
and $B,C,D$ act boundedly. Thus $(AC)|_{Y_x}$ is
quasinilpotent, so
\[
 (BD)|_{Y_x}=I_{Y_x}-(AC)|_{Y_x}
\]
is invertible. Consequently, $B|_{Y_x}$ is invertible.
The sequence
\[
 x_n=(B|_{Y_x})^{-n}x,\qquad n\geq0,
\]
satisfies $Bx_{n+1}=x_n$ and
\[
 \|x_n\|
 \leq\|(B|_{Y_x})^{-1}\|_{\B(Y_x)}^n\|x\|_x.
\]
If $x\ne0$, set
\[
 c=\max\{1,\|(B|_{Y_x})^{-1}\|_{\B(Y_x)}\}
   \max\{1,\|x\|_x/\|x\|\}.
\]
Then $\|x_n\|\leq c^n\|x\|$ for $n\geq0$, and hence
$x\in\K(B)$. The case $x=0$ is immediate. Symmetry proves
\eqref{eq:bezout-cross}.

Clearly,
\[
 \Ran(AB)\subseteq\Ran(A)\cap\Ran(B).
\]
If $x=Au=Bv$, then
\[
 x=(AC+BD)x=AB(Cv+Du),
\]
which proves~\eqref{eq:bezout-range}.

Using~\eqref{eq:bezout-K}--\eqref{eq:bezout-cross},
we have
\[
 L(AB)
 =(\K(A)\cap\K(B))+\Hq(A)+\Hq(B).
\]
If
\[
 x=k_A+h_A=k_B+h_B\in L(A)\cap L(B),
\]
then
\[
 k_A-h_B=k_B-h_A\in\K(A)\cap\K(B).
\]
This proves $L(A)\cap L(B)\subseteq L(AB)$.
The reverse inclusion follows from
 the assertion \eqref{eq:bezout-cross} above. The proof for $G$ is identical,
replacing the analytic cores by ranges and using
 the assertion \eqref{eq:bezout-range} above.

Finally, suppose that $Ax\in G(AB)$. Write
\[
 Ax=ABy+h_A+h_B,
 \qquad h_A\in\Hq(A),\quad h_B\in\Hq(B).
\]
By the B\'ezout identity,
\[
 x=CAx+BDx
 =B(CAy+Dx)+Ch_A+Ch_B\in G(B),
\]
since
\[
 Ch_A\in\K(B)\subseteq\Ran(B),
 \qquad Ch_B\in\Hq(B).
\]
Conversely, if $x=By+h_B\in G(B)$, then
\[
 Ax=ABy+Ah_B\in G(AB).
\]
Thus $A^{-1}(G(AB))=G(B)$, and the other identity
follows by symmetry.
\end{proof}

\begin{proposition}\label{prop:left-regular-bezout}
Let $H$ be a complex Hilbert space.
Then $\LGD(H)$ and $\RGD(H)$ are
regularities in $\B(H)$.
\end{proposition}

\begin{proof}
We first consider $\LGD(H)$.
Let $T\in\B(H)$ and $n\geq1$.
Since
\[
 \Hq(T^n)=\Hq(T),
\]
the subspace $N=\Hq(T)$ is closed
for $T$ if and only if it is closed
for $T^n$. Whenever $N$ is closed,
the operators induced by $T^n$ and
$T$ on $H/N$ satisfy
\[
 \widehat{T^n}=\widehat T^{\,n}.
\]
By Lemma~\ref{lem:left-criterion},
$T\in\LGD(H)$ if and only if $N$
is closed and $\widehat T$ is
bounded below. Since an operator
is bounded below if and only if
any of its positive powers is
bounded below, we obtain
\[
 T\in\LGD(H)
 \quad\Longleftrightarrow\quad
 T^n\in\LGD(H).
\]

Now let $A,B,C,D\in\B(H)$ commute
pairwise and satisfy $AC+BD=I$.

Suppose first that $AB\in\LGD(H)$.
By Lemma~\ref{lem:left-criterion},
$\Hq(AB)$ and $G(AB)$ are closed.
Lemma~\ref{lem:bezout} implies
that $\Hq(A)$ and $\Hq(B)$ are
closed. Moreover, by  the assertion 
\eqref{eq:bezout-G-preimages} above,
\[
 G(A)=B^{-1}G(AB),\qquad
 G(B)=A^{-1}G(AB),
\]
and thus are closed. Hence
$A,B\in\LGD(H)$ by
Lemma~\ref{lem:left-criterion}.

Conversely, suppose that
$A,B\in\LGD(H)$. Then
$\Hq(A)$, $\Hq(B)$, $G(A)$
and $G(B)$ are closed.
Choose a reducing decomposition
\[
 H=M_A\oplus\Hq(A)
\]
on which $A|_{M_A}$ is bounded
below. Since
$
 \K(A|_{\Hq(A)})=\{0\},
$
we have $\K(A)\subseteq M_A$.
By the assertion \eqref{eq:bezout-cross} above,
$\Hq(B)\subseteq M_A$.
Hence
\[
 \Hq(A)\oplus\Hq(B)
\]
is closed. By
Lemma~\ref{lem:bezout},
$\Hq(AB)$ is closed and
\[
 G(AB)=G(A)\cap G(B)
\]
is closed. Another application
of Lemma~\ref{lem:left-criterion}
gives $AB\in\LGD(H)$.

Thus $\LGD(H)$ satisfies both
regularity axioms.

Finally,
\[
 T\in\RGD(H)
 \quad\Longleftrightarrow\quad
 T^*\in\LGD(H).
\]
Since taking adjoints preserves
powers and the commuting B\'ezout
identity, the regularity of
$\RGD(H)$ follows from that
of $\LGD(H)$.
\end{proof}

\begin{theorem}\label{thm:main}
For every complex Hilbert space $H$,
the class $\PF(H)$ of pseudo-Fredholm
operators is a regularity in $\B(H)$.
\end{theorem}

\begin{proof}
The identity operator is semi-regular,
so $\PF(H)$ is nonempty.

We first verify the assertion \eqref{eq:regularity-powers} above.
For each integer $n\geq1$,
\begin{equation}\label{eq:power-cores}
 \K(T^n)=\K(T),\qquad
 \Hq(T^n)=\Hq(T).
\end{equation}
The first equality follows from
\cite[Corollary~2.1]{Tajmouati2014}.

Put $L=\K(T)+\Hq(T)$, and suppose
that $L$ is closed. By
\eqref{eq:core-facts},
$T^{-1}L=L$; therefore the operator
$\widehat T$ induced on $H/L$ is
injective. Since
$\K(T)\subseteq\Ran(T^n)$,
\[
 \Ran(\widehat T^{\,n})
 =
 \frac{\Ran(T^n)+L}{L}
 =
 \frac{\Ran(T^n)+\Hq(T)}{L}.
\]
As $L$ is closed,
$\Ran(T^n)+\Hq(T)$ is closed
if and only if
$\Ran(\widehat T^{\,n})$
is closed. Since $\widehat T$
is injective, this is equivalent
to $\widehat T^{\,n}$ being
bounded below, which is equivalent
to $\widehat T$ being bounded below.
Consequently,
\[
 \Ran(T^n)+\Hq(T)\text{ is closed}
 \quad\Longleftrightarrow\quad
 \Ran(T)+\Hq(T)\text{ is closed}.
\]
Together with the assertion \eqref{eq:power-cores} above
and Proposition~\ref{prop:criterion},
this proves
\[
 T^n\in\PF(H)
 \quad\Longleftrightarrow\quad
 T\in\PF(H).
\]

Now let $A,B,C,D\in\B(H)$ commute
pairwise and satisfy $AC+BD=I$.
Suppose that $A,B\in\PF(H)$.
The three subspaces
\[
 \K(A),\quad L(A),\quad G(A)
\]
and the corresponding subspaces
for $B$ are closed by the assertion 
\eqref{eq:hilbert-criterion}.
Lemma~\ref{lem:bezout} gives
\[
 \K(AB)=\K(A)\cap\K(B), \,\,
 L(AB)=L(A)\cap L(B),\,\,
 G(AB)=G(A)\cap G(B).
\]
All three subspaces are closed,
and hence $AB\in\PF(H)$.

Conversely, suppose that $T=AB$
is pseudo-Fredholm and put
$E=\K(T)$. By
Lemma~\ref{lem:quotient},
$E$ is closed, $T(E)=E$, and
the operator $\widehat T$
induced on $H/E$ is left
generalized Drazin invertible.
Every operator commuting with
$T$ leaves $\K(T)$ invariant.
Thus $A,B,C,D$ induce commuting
operators on $H/E$ satisfying
\[
 \widehat A\widehat C+
 \widehat B\widehat D=I,
 \qquad
 \widehat A\widehat B=\widehat T.
\]
Since $H/E$ is a Hilbert space,
Proposition~\ref{prop:left-regular-bezout}
shows that $\widehat A$ and
$\widehat B$ are left generalized
Drazin invertible. Moreover,
\[
 E=AB(E)\subseteq A(E)\subseteq E,
\]
\[
 E=AB(E)\subseteq B(E)\subseteq E.
\]
Thus $A(E)=B(E)=E$.
Lemma~\ref{lem:quotient}, applied
to $A$ and $B$, yields
$A,B\in\PF(H)$.
\end{proof}

\begin{corollary}\label{cor:spectral-mapping}
Let $H$ be a complex Hilbert space.
For $T\in\B(H)$, put
\[
\begin{aligned}
 \sigma_{gK}(T)
 &=\{\lambda\in\mathbb C:
 T-\lambda I\notin\PF(H)\},\\
 \sigma_{lgD}(T)
 &=\{\lambda\in\mathbb C:
 T-\lambda I\notin\LGD(H)\},\\
 \sigma_{rgD}(T)
 &=\{\lambda\in\mathbb C:
 T-\lambda I\notin\RGD(H)\},\\
 \sigma_{gD}(T)
 &=\{\lambda\in\mathbb C:
 T-\lambda I\notin\GD(H)\}.
\end{aligned}
\]
Then
\[
 \sigma_{\bullet}(f(T))
 =f\bigl(\sigma_{\bullet}(T)\bigr),
 \qquad
 \bullet\in\{gK,lgD,rgD,gD\},
\]
for every function $f$ holomorphic
on a neighbourhood of $\sigma(T)$
and nonconstant on each connected
component of that neighbourhood.
\end{corollary}

\begin{proof}
By Theorem~\ref{thm:main} and
Proposition~\ref{prop:left-regular-bezout},
the classes $\PF(H)$, $\LGD(H)$
and $\RGD(H)$ are regularities.

Since $\GD(H)=\LGD(H)\cap\RGD(H)$
(see \cite[Proposition~2.18(iii)--(v)]{AznayOuahabZariouh2023}),
$\GD(H)$ is also a regularity. The conclusion
follows from the spectral mapping theorem for regularities
\cite[Theorem~1.4]{KM1996}.
\end{proof}

\section{Banach-space counterexample and spectral regularities}

The power axiom fails on some separable Banach spaces, as the
following adaptation of
\cite[Example~5.3 and the proof of Theorem~5.5]{AznayChabbabi2026}
shows.

\begin{theorem}\label{thm:counterexample}
There exist a separable complex Banach
space $X$ and an operator $T\in\B(X)$
such that $T^2$ is right generalized
Drazin invertible, whereas $T$ is
not pseudo-Fredholm. In particular,
neither $\PF(X)$ nor the class of
right generalized Drazin invertible
operators need be a regularity in
$\B(X)$.
\end{theorem}

\begin{proof}
Choose a bounded surjection $q:\ell^1\to\ell^2$,
whose existence follows from
\cite[Theorem~6.1]{Carothers2004}.

Put $E=\ell^1\oplus_1\ell^2$
and $F=\ell^2$. For
$a=(a_j)\in\ell^1$ write
\[
 a_o=(a_{2j-1})_{j\geq1},
 \qquad
 a_e=(a_{2j})_{j\geq1}.
\]
Define $S\in\B(E)$ and
$J\in\B(F,E)$ by
\[
 S(a,b)=(a_o,q(a_e)),
 \qquad Jy=(0,y).
\]
The operator $S$ is surjective:
to solve $S(a,b)=(u,v)$, choose
$d\in\ell^1$ such that $q(d)=v$
and interleave $u$ and $d$ as
the odd and even coordinates
of $a$. Moreover, $SJ=0$.

On $X=E\oplus_1F$ define
\[
 T=
 \begin{pmatrix}
 S&J\\
 0&0
 \end{pmatrix}.
\]
Then
\[
 T^2=
 \begin{pmatrix}
 S^2&0\\
 0&0
 \end{pmatrix}.
\]
Since $S^2$ is surjective,
$T^2$ is the direct sum of
a surjective operator and a
nilpotent operator. In particular,
$T^2\in\RGD(X)\subseteq\PF(X)$.

Since $S$ is surjective, $\K(S)=E$;
see \cite[p.~3523]{GonzalezMbekhtaOudghiri2008}.
On the other hand,
\[
 \Ran(T)\subseteq E\oplus\{0\}.
\]
Consequently,
\begin{equation}\label{eq:example-core}
 \K(T)=E\oplus\{0\}.
\end{equation}

Suppose that $T$ were
pseudo-Fredholm, with a generalized
Kato decomposition $X=M\oplus N$.
Since $\K(T|_N)=\{0\}$,
\eqref{eq:example-core} implies
$E\oplus\{0\}\subseteq M$.

If $(x,y)\in M$ with $y\ne0$,
choose $u\in E$ such that
\[
 S(x+u)+Jy=0,
\]
which is possible because $S$
is surjective. Then
$(x+u,y)\in M\cap\Ker(T)$.
As
\[
 \Ran(T|_M)\subseteq E\oplus\{0\},
\]
this vector does not belong
to $\Ran(T|_M)$, contrary to
the semi-regularity of $T|_M$.
Thus
\[
 M=E\oplus\{0\}.
\]

Every closed complement $N$
of $M$ is the graph of a
bounded operator $V:F\to E$.
Since $N$ is $T$-invariant
and $T(N)\subseteq M$, we
have $T(N)=\{0\}$, that is,
\begin{equation}\label{eq:example-rightinverse}
 SV=-J.
\end{equation}
Write $V=(V_1,V_2)$, where
$V_1:F\to\ell^1$, and let
$P_e:\ell^1\to\ell^1$ extract
the even coordinates. Taking
the second coordinate in
\eqref{eq:example-rightinverse}
gives
\[
 qP_eV_1=-I_{\ell^2}.
\]
Let $(e_j)$ be the canonical
basis of $\ell^2$. Since
$e_j\to0$ weakly, the sequence
$(P_eV_1e_j)$ tends weakly to
zero in $\ell^1$. The Schur
property of $\ell^1$ implies
\[
 \|P_eV_1e_j\|_1\to0,
\]
contradicting
\[
 \|qP_eV_1e_j\|_2=1.
\]
Hence $T$ is not
pseudo-Fredholm.
\end{proof}

For $T\in\B(X)$, let $\sigma_{ap}(T)$ and $\sigma_{su}(T)$
denote the approximate point and surjectivity spectra,
respectively.

\begin{remark}\label{rem:surjectivity-isolated}
For the operator in Theorem~\ref{thm:counterexample},
$0\in\operatorname{iso}\sigma_{su}(T)$ although $T$ is not right
generalized Drazin invertible. Indeed, surjectivity is open in
$\B(E)$, so $S-\lambda I$ is surjective for all sufficiently
small $|\lambda|$. For such $\lambda\ne0$,
\[
 T-\lambda I=\begin{pmatrix}S-\lambda I&J\\0&-\lambda I\end{pmatrix}
\]
is surjective, whereas $T$ is not surjective.

The same failure occurs on Hilbert spaces. The bilateral
weighted shift $B$ of \cite[Proposition~2.1]{AznayChabbabi2026}
satisfies $\K(B)+\Hq(B)=\ell^2(\mathbb Z)$ and is not
pseudo-Fredholm. Since $B$ is not surjective,
\cite[Theorem~5]{GonzalezMbekhtaOudghiri2008} gives
$0\in\operatorname{iso}\sigma_{su}(B)$. Thus the equivalence
between right generalized Drazin invertibility and
$0\notin\operatorname{acc}\sigma_{su}(T)$ asserted in
\cite[Theorem~3.8]{BenharratHocineMessirdi2019}
fails even on Hilbert spaces.
\end{remark}

The classes of bounded-below and surjective operators are
regularities; see \cite[Example~(5)]{KM1996}.

\begin{proposition}\label{prop:acc-one-sided}
For every complex Banach space $X$, the two classes
\[
\begin{aligned}
 \mathcal R_{ap}^{acc}(X)
 &=\{T\in\B(X):0\notin\operatorname{acc}\sigma_{ap}(T)\},\\
 \mathcal R_{su}^{acc}(X)
 &=\{T\in\B(X):0\notin\operatorname{acc}\sigma_{su}(T)\}
\end{aligned}
\]
are regularities in $\B(X)$. Their associated spectra are
$\operatorname{acc}\sigma_{ap}(T)$ and
$\operatorname{acc}\sigma_{su}(T)$, respectively.
\end{proposition}

\begin{proof}
We prove the assertion for $\mathcal R_{ap}^{acc}(X)$ first.
This class contains $I$. By the polynomial spectral mapping theorem
\cite[Theorem~1]{ChoiDavis1974}, for every $n\geq1$,
\[
 \sigma_{ap}(T^n)=\{\lambda^n:\lambda\in\sigma_{ap}(T)\}.
\]
Since taking positive powers preserves accumulation at zero, the
power axiom follows.

Let $A,B,C,D\in\B(X)$ commute pairwise and satisfy $AC+BD=I$.
Put
\[
 \Omega=\sigma_{ap}(A,B)
 =\left\{(\alpha,\beta)\in\mathbb C^2:
   \inf_{\|x\|=1}\bigl(\|(A-\alpha I)x\|
                         +\|(B-\beta I)x\|\bigr)=0\right\},
\]
the joint approximate point spectrum of $(A,B)$. By
\cite[Theorem~1]{ChoiDavis1974}, $\Omega$ is compact,
its coordinate projections are $\sigma_{ap}(A)$ and
$\sigma_{ap}(B)$, and
\[
 \sigma_{ap}(AB)=\{\alpha\beta:(\alpha,\beta)\in\Omega\}.
\]
For $(\alpha,\beta)\in\Omega$, choose unit vectors $x_k$ such that
$(A-\alpha I)x_k\to0$ and $(B-\beta I)x_k\to0$. Since
\[
 x_k=C(A-\alpha I)x_k+D(B-\beta I)x_k
       +\alpha Cx_k+\beta Dx_k,
\]
we obtain
\[
 1\leq |\alpha|\,\|C\|+|\beta|\,\|D\|.
\]

Suppose $0\in\operatorname{acc}\sigma_{ap}(A)$. Choose
$\alpha_k\in\sigma_{ap}(A)\setminus\{0\}$ with $\alpha_k\to0$.
By the projection property, there are $\beta_k$ such that
$(\alpha_k,\beta_k)\in\Omega$. The preceding inequality yields
$\beta_k\ne0$ for all sufficiently large $k$. Moreover,
$|\beta_k|\leq\|B\|$, so $0\ne\alpha_k\beta_k\to0$ and
$0\in\operatorname{acc}\sigma_{ap}(AB)$. The argument for $B$
is symmetric.

Conversely, suppose $0\in\operatorname{acc}\sigma_{ap}(AB)$.
Choose $(\alpha_k,\beta_k)\in\Omega$ with
$\alpha_k\beta_k\ne0$ and $\alpha_k\beta_k\to0$.
By compactness, after passing to a subsequence, we have
$(\alpha_k,\beta_k)\to(\alpha,\beta)\in\Omega$.
Then $\alpha\beta=0$, while the preceding inequality gives
$(\alpha,\beta)\ne(0,0)$. If $\alpha=0$, then
$\alpha_k\ne0$ and $\alpha_k\to0$, so
$0\in\operatorname{acc}\sigma_{ap}(A)$; the case $\beta=0$
is analogous. Consequently,
\[
0\in\operatorname{acc}\sigma_{ap}(AB)\iff 0\in\operatorname{acc}\sigma_{ap}(A)\ \text{or}\ 0\in\operatorname{acc}\sigma_{ap}(B).
\]

Finally, $\sigma_{su}(T)=\sigma_{ap}(T^*)$ on every Banach space.
As adjoints preserve powers and the commuting B\'ezout identity,
the two axioms for $\mathcal R_{su}^{acc}(X)$ follow from those for
$\mathcal R_{ap}^{acc}(X^*)$. The associated spectra are identified
by $\sigma_{\bullet}(T-\lambda I)=\sigma_{\bullet}(T)-\lambda$,
for $\bullet\in\{ap,su\}$.
\end{proof}

We also obtain another proof of Lubansky's regularity theorem
\cite{Lubansky2007} for $\B(X)$.

\begin{corollary}\label{cor:banach-gd-regularity}
For every complex Banach space $X$,
\[
 \GD(X)=\mathcal R_{ap}^{acc}(X)
       \cap\mathcal R_{su}^{acc}(X).
\]
In particular, $\GD(X)$ is a regularity in $\B(X)$.
Consequently, for every function $f$ holomorphic on a
neighbourhood of $\sigma(T)$ and nonconstant on each
connected component of that neighbourhood,
\[
 \operatorname{acc}\sigma_{\bullet}(f(T))
 =f\bigl(\operatorname{acc}\sigma_{\bullet}(T)\bigr),
 \qquad \sigma_{\bullet}\in\{\sigma,\sigma_{ap},\sigma_{su}\}.
\]
\end{corollary}

\begin{proof}
Since $\sigma(T)=\sigma_{ap}(T)\cup\sigma_{su}(T)$,
\[
 \operatorname{acc}\sigma(T)
 =\operatorname{acc}\sigma_{ap}(T)
  \cup\operatorname{acc}\sigma_{su}(T).
\]
The spectral characterization of generalized Drazin
invertibility, $T\in\GD(X)$ if and only if
$0\notin\operatorname{acc}\sigma(T)$, gives the equality.
The intersection of two regularities is a regularity.
The spectral mapping identities follow from
\cite[Theorem~1.4]{KM1996}.
\end{proof}

\begin{remark}
The regularity $\mathcal R_{su}^{acc}(X)$ must not be confused
with the class $\RGD(X)$. They differ for the operator of
Theorem~\ref{thm:counterexample}, by
Remark~\ref{rem:surjectivity-isolated}.
\end{remark}

\end{document}